\documentclass[a4paper, 11pt]{article}

\usepackage[T1]{fontenc}
\usepackage[utf8]{inputenc}
\usepackage{lmodern, microtype}
\usepackage{latexsym, textcomp}
\usepackage{amsmath, amssymb, amsthm}
\usepackage{array, setspace, verbatim}
\usepackage{xargs, xifthen}

\usepackage[english]{babel}
\usepackage[english]{varioref}

\usepackage[pdftex]{hyperref}
\usepackage{pdfpages}

\newcommand*{\nidt}{\noindent}
\newcommand*{\dst}{\displaystyle}
\newcommand*{\ov}[1]{\overline{#1}}
\newcommand*{\wh}[1]{\widehat{#1}}

\newcommand*{\esp}{\vspace*{0.5ex}}
\newcommand*{\espf}{\vspace*{1ex}}

\newcommand*{\ii}{\infty}
\newcommand*{\del}{\nabla}
\renewcommand*{\phi}{\varphi}

\newcommand*{\C}{\mathbb{C}}
\renewcommand*{\H}{\mathbb{H}}
\newcommand*{\Mm}{\mathcal{M}}
\newcommand*{\R}{\mathbb{R}}
\renewcommand*{\S}{\mathbb{S}}

\newcommand*{\Nil}{\text{Nil}}

\newcommand*{\lp}{\left(}
\newcommand*{\rp}{\right)}
\newcommand*{\lc}{\left[}
\newcommand*{\rc}{\right]}
\newcommand*{\lac}{\left\{}

\newcommand*{\lan}{\langle}
\newcommand*{\ran}{\rangle}

\newcommand*{\rr}{\right.}

\newcommand*{\into}{\rightarrow}

\renewcommand*{\matrix}[2]{\lp\begin{array}{@{}*{#1}{c}@{}}#2\end{array}\rp}
\newcommand*{\der}[2][]{\frac{\partial#1}{\partial#2}}

\newtheorem{theorem}{Theorem}[section]

\newtheorem{e-proposition}[theorem]{Proposition}

\newtheorem{e-definition}[theorem]{Definition\rm}

\hypersetup{%
pdfauthor= {S{\'e}bastien Cartier},%
pdftitle= {Sym-Bobenko formula for minimal surfaces in Heisenberg space},%
pdfcreator= {PDFLaTeX},%
pdfproducer= {PDFLaTeX}%
}

\title{Sym-Bobenko formula for minimal\\surfaces in Heisenberg space\footnote{This work is part of the author's Ph.~D. thesis \cite{Ca}.}}
\author{S{\'e}bastien Cartier}

\begin{document}

\maketitle

\begin{abstract} \itshape

We give an immersion formula, the \emph{Sym-Bobenko formula}, for minimal surfaces in the $3$-dimensional Heisenberg space. Such a formula can be used to give a generalized Weierstrass type representation and construct explicit examples of minimal surfaces.

\end{abstract}

\nidt \textit{Mathematics Subject Classification:} \emph{Primary 53A10, Secondary 53C42}.

\section{Introduction}

A Sym-Bobenko formula is the expression of an immersion in terms of a one-parameter family of moving frames, called the \emph{extended frame}. This idea was first used by A.~Sym \cite{Sy} in the case of surfaces with negative constant (Gauss) curvature in euclidean space. A.~I.~Bobenko applied the method to numerous cases \cite{Bo2} \cite{Bo} \cite{Bo3}, including constant mean curvature (\emph{CMC} for short) surfaces in space forms ---~euclidean $3$-space, $3$-sphere and hyperbolic $3$-space~--- and T.~Taniguchi applied it to CMC spacelike surfaces in Minkowski $3$-space \cite{Ta}. In the mean time, the work of J.~Dorfmeister, F.~Pedit and H.~Wu \cite{DoPeWu} and D.~Brander, W.~Rossman and N.~Schmitt \cite{BrRoSc} show that Sym-Bobenko formul{\ae} can be seen as generalized Weierstrass type representations for CMC surfaces, extended frames coming from holomorphic data.

In Heisenberg $3$-space, the classical method does not apply, since the isometry group is of dimension only $4$ ---~contrary to the ones of space forms that are $6$-dimensional~--- and does not act transitively on orthonormal frames; there are ``not enough'' isometries to define a moving frame. We show that nevertheless, for minimal immersions a Sym-Bobenko formula can be established using an \textit{ad-hoc} matrix-valued map.

In \cite{DoInKo2}, J.~F.~Dorfmeister, J.~Inoguchi and S.~Kobyashi link this formula with pairs of meromorphic and anti-meromorphic $1$-forms, which they call \emph{pairs of normalized potentials}, in a way to get a generalized Weierstrass type representation for minimal surfaces.

\section{Surfaces in Heisenberg space}

We see the $3$-dimensional Heisenberg space $\Nil_3$ as $\R^3$, with generic coordinates $(x_1, x_2, x_3)$, endowed with the following riemannian metric:
\[
\lan \cdot, \cdot \ran= dx_1^2+ dx_2^2+ \lp \frac{1}{2} (x_2 dx_1- x_1 dx_2)+ dx_3 \rp^2.
\]
We call \emph{canonical frame} the orthonormal frame $(E_1, E_2, E_3)$ defined by:
\[
E_1= \der{x_1}- \frac{x_2}{2} \der{x_3}, \quad E_2= \der{x_2}+ \frac{x_1}{2} \der{x_3} \quad \text{and} \quad E_3= \der{x_3},
\]
and the Levi-Civita connection $\del$ writes:
\[
\begin{array}{l@{\quad}l@{\quad}l}
\del_{E_1} E_1= 0                & \del_{E_2} E_1= -\frac{1}{2} E_3 & \del_{E_3} E_1= -\frac{1}{2} E_2 \espf \\
\del_{E_1} E_2= \frac{1}{2} E_3  & \del_{E_2} E_2= 0                & \del_{E_3} E_2= \frac{1}{2} E_1  \espf \\
\del_{E_1} E_3= -\frac{1}{2} E_2 & \del_{E_2} E_3= \frac{1}{2} E_1  & \del_{E_3} E_3= 0.
\end{array}
\]
Note that the vector field $E_3$ is a Killing field and that the projection $\pi: (x_1, x_2, x_3) \in \Nil_3 \mapsto (x_1, x_2) \in \R^2$ on the first two coordinates is a Riemannian submersion. From now on, we identify $\R^2$ with $\C$.

We may also write $\Nil_3$ as a subset of $\Mm_2(\C)$. Consider the matrices:
\[
\sigma_0= \matrix{2}{i & 0 \\ 0 & -i}, \quad \sigma_1= \matrix{2}{0 & 1 \\ 1 & 0}, \quad \sigma_2= \matrix{2}{0 & i \\ -i & 0} \quad \text{and} \quad \sigma_3= \matrix{2}{1 & 0 \\ 0 & 1}.
\]
The identification is the following:
\begin{multline} \label{eq:matrixmodel}
(x_1, x_2, x_3) \in \Nil_3 \longleftrightarrow x_1 \sigma_1+ x_2 \sigma_2+ x_3 \sigma_3 \\
= \matrix{2}{x_3 & x_1+ ix_2 \\ x_1- ix_2 & x_3} \in \Mm_2(\C).
\end{multline}
Note that this identification is purely formal and does not involve any manifold related structure.

\medskip

Let $\Sigma$ be a simply connected Riemann surface and $z$ be a conformal parameter on $\Sigma$. A conformal immersion is denoted $f: \Sigma \into \Nil_3$ with unit normal $N$ and conformal factor $\rho: \Sigma \into (0, +\ii)$ meaning:
\begin{gather*}
\lan f_z, f_z \ran= \lan f_{\bar{z}}, f_{\bar{z}} \ran= 0, \quad \lan f_z, f_{\bar{z}} \ran= \frac{\rho}{2}, \\
\lan f_z, N \ran= \lan f_{\bar{z}}, N \ran= 0 \quad \text{and} \quad \lan N, N \ran= 1.
\end{gather*}
Consider also $\phi= \lan N, E_3 \ran: \Sigma \into (-1, 1)$ denote the angle function of $N$, $A= \lan f_z, E_3 \ran: \Sigma \into \C$ the vertical part of $f_z$ and $pdz^2= \lan \del_{f_z} f_z, N \ran dz^2$ the Hopf differential of $f$.

The Abresch-Rosenberg differential expresses $Qdz^2= (ip+ A^2)dz^2$, and a necessary and sufficient condition for $f$ to be minimal is $\del_{f_z} f_{\bar{z}}= 0$.

We also decompose $f$ into $f= (F, h)$ with $F= \pi \circ f: \Sigma \into \C$ the horizontal projection of $f$ and $h: \Sigma \into \R$ its height function. We can express $A$ in terms of $F$ and $h$:
\begin{equation} \label{eq:verticalpart}
A= h_z- \frac{i}{4} \lp F\ov{F}_z- \ov{F} F_z \rp.
\end{equation}
In the matrix model~\eqref{eq:matrixmodel} of $\Nil_3$, the map $F$ is given by the non-diagonal coefficients ---~precisely the $(1, 2)$-coefficient~--- and $h$ by the diagonal ones.

\medskip

The intuitive idea behind Sym-Bobenko formul{\ae} in space forms is that, up to ambient isometries, the unit normal ---~or Gauss map~--- would locally determine the immersion up to ambient isometries. In $\Nil_3$ such a map is defined as follows; see \cite{Da3} for details. Since $\Nil_3$ is a Lie group, the map $f^{-1} N$ takes values in the unit sphere $\S^2$ of the Lie algebra. Moreover, for a local study, we can suppose $\phi> 0$ so that the values of $f^{-1} N$ are actually in the northern hemisphere of $\S^2$. If $s$ denotes the stereographic projection centered at the South Pole, we call \emph{Gauss map} of an immersion $f$ the map $g= s \circ (f^{-1} N)$ with values in the unit disk. Actually, endowing the unit disk with the Poincar{\'e} metric, we see the Gauss map $g$ as a map with values into the hyperbolic disk $\H^2$.

\medskip

We use the following criterion to show that a conformal immersion in Heisenberg space is minimal:

\medskip

\begin{e-proposition}[Daniel \cite{Da3}] \label{prop:minimalimmersion}

A conformal immersion $f= (F, h): \Sigma \into \Nil_3$ is minimal if and only if:
\[
F_{z\bar{z}}= \frac{i}{2} \lp \ov{A} F_z+ A F_{\bar{z}} \rp \quad \text{and} \quad A_{\bar{z}}+ \ov{A}_z= 0.
\]
Furthermore, when $f$ is minimal its Gauss map $g: \Sigma \into \H^2$ is harmonic.

\end{e-proposition}

\section{The Sym-Bobenko formula}

Consider the family $(\Psi_t)_{t \in \R}$ of matrix fields over $\Sigma$ which are solutions of the system:
\[
\lac \begin{array}{l}
\dst \Psi_t^{-1} d\Psi_t= \frac{1}{4} \matrix{2}{(\log \rho_0)_z & i\sqrt{\rho_0} \esp \\ \dst -\frac{4iQ_0}{\sqrt{\rho_0}} e^{2it} & -(\log \rho_0)_z \espf} dz \espf \\
\dst \hspace{12em}+ \frac{1}{4} \matrix{2}{-(\log \rho_0)_{\bar{z}} & \dst \frac{4i\ov{Q_0}}{\sqrt{\rho_0}} e^{-2it} \esp \\ -i\sqrt{\rho_0} & (\log \rho_0)_{\bar{z}} \espf} d\bar{z} \\
\Psi_t(z= 0)= \sigma_3
\end{array} \rr,
\]
where $\rho_0: \Sigma \into (0, +\ii)$ and $Q_0: \Sigma \into \C$ are smooth. Such a family $(\Psi_t)$ exists if and only if:
\[
(\log \rho_0)_{z\bar{z}}= \frac{\rho_0}{8}- \frac{2|Q_0|^2}{\rho_0} \quad \text{and} \quad (Q_0)_{\bar{z}}= 0.
\]

\medskip

\begin{theorem}[Sym-Bobenko formula]

Using the matrix model~\eqref{eq:matrixmodel}, define the map $f_t: \Sigma \into \Nil_3$ for any $t \in \R$ as:
\begin{equation} \label{eq:symbobenkoformula}
f_t= -\frac{1}{2} \lp \sigma_0 \der[\wh{f}_t]{t} \rp^d+ \lp \wh{f}_t \rp^{nd} \quad \text{with} \quad \wh{f}_t= -2\der[\Psi_t]{t} \Psi_t^{-1}+ 2 \Psi_t \sigma_0 \Psi_t^{-1},
\end{equation}
where the superscripts $\cdot^d$ and $\cdot^{nd}$ denote respectively the diagonal an non-diagonal terms. Then $f_t$ is a conformal minimal immersion in Heisenberg space and the family $(f_t)$ is the so-called \emph{associated family}.
\end{theorem}

\medskip

\nidt \textit{Proof}. Fix $t \in \R$. From Equation~\eqref{eq:symbobenkoformula}, we get:
\begin{equation} \label{eq:matrixlink}
\begin{array}{@{}c@{}}
\dst (f_t)^{nd}= (\wh{f}_t)^{nd}, \quad \big{(} (f_t)_z \big{)}^d= \frac{1}{2} \lp \sigma_0 \lc (\wh{f}_t)_z, \wh{f}_t \rc \rp^d+ 2\lp \sigma_0 (\wh{f}_t)_z \rp^d \espf \\
\dst \text{and} \quad (\wh{f}_t)_{z\bar{z}}= \frac{i}{4} \lc (\wh{f}_t)_z, (\wh{f}_t)_{\bar{z}} \rc,
\end{array}
\end{equation}
where $[\cdot, \cdot]$ denotes the commutator. From the first equation in~\eqref{eq:matrixlink}, we have that matrices $f_t$ and $\wh{f}_t$ write:
\[
f_t= \matrix{2}{h & F \\ \ov{F} & h} \quad \text{and} \quad \wh{f}_t= \matrix{2}{i\wh{h} & F \\ \ov{F} & -i\wh{h}},
\]
with $F: \Sigma \into \C$ and $h, \wh{h}: \Sigma \into \R$ smooth. We show that $F$ and $h$ verify the conditions of Proposition~\ref{prop:minimalimmersion}. Using Equation~\eqref{eq:verticalpart} and the second identity in~\eqref{eq:matrixlink}, we deduce $A= i\wh{h}_z$ and since $\wh{h}$ is real-valued, we obtain $A_{\bar{z}}+ \ov{A}_z= 0$. Finally, the $(1, 2)$-coefficient of the third equation in~\eqref{eq:matrixlink} verify:
\[
F_{z\bar{z}}= \frac{1}{2} \lp \wh{h}_{\bar{z}} F_z- \wh{h}_z F_{\bar{z}} \rp= \frac{i}{2} \lp \ov{A} F_z+ A F_{\bar{z}} \rp,
\]
which concludes the proof. \hfill$\square$

\bibliographystyle{amsplain}
\providecommand{\MR}{\relax\ifhmode\unskip\space\fi MR }
\providecommand{\MRhref}[2]{%
  \href{http://www.ams.org/mathscinet-getitem?mr=#1}{#2}
}
\providecommand{\href}[2]{#2}

\vfill

\nidt Universit{\'e} Paris-Est, Laboratoire d'Analyse et de Math{\'e}matiques Appliqu{\'e}es (UMR 8050), UPEC, UPEMLV, CNRS, F-94010, Cr{\'e}teil, France \\
\textit{E-mail address:} \verb+sebastien.cartier@u-pec.fr+

\end{document}